\documentclass[11pt]{article}
\usepackage{makeidx}
\usepackage{latexsym}
\usepackage{graphicx}
\usepackage[portuguese,english]{babel}
\usepackage[latin1]{inputenc}
\usepackage[T1]{fontenc}
\usepackage[active]{srcltx}
\usepackage{pstricks,pst-node,pst-text,pst-3d,pst-plot}
\makeindex
\newtheorem{teo}{Theorem}

\newtheorem{lema}{Lemma}

\begin{document}

\centerline{\textbf{How to find finite topological spaces}}

\centerline{\textbf{with given quotient-spaces}}

\vspace{0.2cm}

\centerline{J. M. S. Sim\~{o}es-Pereira}

\centerline{(orcid: 0000.0001.8930.5932)}

\vspace{0.1cm}

\centerline{Department of Mathematics, University of Coimbra, Coimbra, Portugal}
\centerline{E-mail: siper@mat.uc.pt}

\vspace{0.1cm}

\begin{flushright}\emph{To my dear Tumim, in Memoriam}\end{flushright}

\vspace{0.1cm}

\emph{\textbf{Abstract:}} \emph{Our main problem is to find a finite topological space to within homeomorphism, given (also to within homeomorphism) the quotient spaces obtained by identifying one point of the space with each one of the other points. Initially, our aim was to reconstruct a topological space from its quotient-spaces; but a reconstruction is not always possible either in the sense that several non-homeomorphic topological spaces yield the same quotient-spaces, or in the sense that no topological space yields an arbitrarily given family of quotient-spaces. We present an algorithm that detects, and deals with, all these situations.}

\vspace{0.1cm}

\textbf{\emph{Key words:}} quotient-spaces, finite topologies, reconstruction procedures.

\vspace{0.1cm}

\section{Introduction}

Let $T$ be a set, whose elements we are going to call points, let $x,y\in T$ and $T(x,y)=\{\{z\}| z\in T-\{x,y\}\}\cup\{\{x,y\}\}$.
Consider the map $f:T\rightarrow T(x,y)$ (sometimes called the natural map) defined by setting $f(z)=\{x,y\}$ when $z\in \{x,y\}$ and $f(z)=\{z\}$ when $z\in T-\{x,y\}$. In colloquial language, we say that $T(x,y)$ was obtained from $T$ by an identification of $x$ and $y$.

If $\cal T$ is a topology on $T$, the pair $(T,\cal T)$ is said to be a topological space; and the topology obtained by identifying $x$ and $y$, denoted ${\cal T}_{x,y}$, is defined as the topology on  $T(x,y)$  for which a set $S\in {\cal T}_{x,y}$ is open if and only if $f^{-1}(S)$ is an open set in $\cal T$. We refer to this operation as a \emph{topological identification} and  $(T(x,y), {\cal T}_{x,y}) $ is the topological \emph{quotient-space} obtained by this identification.

We deal with the following

\vspace{0.2cm}

\emph{Problem:} Let $(Q_{1},{\cal T}_{1}), ..., (Q_{n-1},{\cal T}_{n-1})$ be $n-1$ given topological spaces,           each one of them with $n-1$ points. Find a topological space $(T, \cal T)$ with $T=\{1, ..., n\}$ such that the given spaces are homeomorphic to those obtained by topological identification of one point of $T$, say $n$ without loss of generality, with each one of the other points of $T$.

\vspace{0.1cm}

This problem may have several non-homeomorphic solutions.

\vspace{0.1cm}

\emph{Example 1:} Let $T=\{1,2,3,4,5\}$ and $n=5$. Let ${\cal T}_{1}$ have  $\emptyset $, $\{1,2\}$ and $T$ as open sets and ${\cal T}_{2}$ have  $\emptyset $, $\{1,2,5\}$ and $T$ as open sets. Both yield two quotient-spaces with $\emptyset$, $\{u,v\}$ and $\{u,v,x,z\}$ as open sets plus two other quotient-spaces with $\emptyset$ and $\{u,v,x,z\}$ as open sets. To each one of these topologies we may add the set $T-\{5\}$ as open set: the quotient-spaces remain the same.

\vspace{0.1cm}

Another case where several topologies yield the same quotient-spaces:

\vspace{0.1cm}

\emph{Example 2:} With $T=\{1,...,7\}$ and $n=7$, the following three topologies ${\cal T}$ with $\emptyset $, $\{1,2,3\}$, $\{4,5,6\}$, $T-\{7\}$ and $T$ as open sets, ${\cal T'}$ with $\emptyset $, $\{1,2,3\}$, $\{1,2,3,7\}$ and $T$ as open sets, and ${\cal T''}$ with $\emptyset $, $\{1,2,3\}$, $\{1,2,3,7\}$, $T-\{7\}$ and $T$ as open sets, yield 6 quotient-spaces all of them with the same topology, namely, $\emptyset $, $\{u,v,w\}$ and $\{u,v,w,x,y,z\}$ as open sets.

\vspace{0.1cm}

Note that we may have situations where no solution exists. For instance, it is a consequence of Theorem 2, that, if, among the quotient-spaces, one of them has no singleton as open set and the remaining ones have three singletons each, then no solution exists.

\vspace{0.1cm}

Initially, we were assuming the  $n-1$ given spaces as quotient-spaces of a finite topological space on $n$ points: our problem was to find it, that means, to reconstruct it. Meanwhile, we designed an algorithmic procedure which does achieve this goal or shows us that no such space exists.

\vspace{0.2cm}

The reconstruction of finite structures when some substructures or some related structures are given is a topic which has captured the attention of several authors. In Graph Theory, the Ulam conjecture is the most famous case (see, for instance, \cite{Harary}); but see also several papers like  \cite{Bondy}, \cite{Bondy2} and \cite{Alice}.
Concerning finite topologies, in \cite{Das}, a topological $n-$space supposed to be connected, $T_{0}$ and $T_{5}$, was reconstructed (to within homeomorphism),  given (also to within homeomorphism) the subspaces induced on its $(n-1)-$subsets. The statement of our problem is similar but, obviously, the problem is not the same! We have already dealt with a much weaker version in \cite{papermeu}.

\vspace{0.1cm}

For Topology, we cite classical texts like \cite{Bushaw} or \cite{Kowalsky} or more recent ones like \cite{Krantz},  \cite{meu} and \cite{Sutherland}. For finite topological spaces, Stong \cite{Stong} was a pioneer, but the topic did not catch the attention of many followers.

\vspace{0.1cm}

As applications for this type of problem, I believe they may be found in the Biomedical Sciences.
Think of manipulating DNA segments of a gene, like fusing or pasting them and how to recover the gene from such manipulated sequences.
These procedures evoke identifications and other operations in finite topological spaces. We don't pretend that immediate, direct applications already exist, but neither are we skeptical about their future existence.

\section{Notation}

As already announced, with no loss of generality we always suppose that $y=n$. While dealing with our \emph{Problem},
we are given a family $F$ of $n-1$ topological spaces (each one defined on a set with $n-1$ points) which we denote by  $Q_{1}$, ..., $Q_{n-1}$, with subscripts; a superscript, say $Q^{i}$, means that this space is known to be the quotient-space obtained from $(T,\cal T)$ by topological identification of points $n$ and $i$: hence, superscripts and subscripts have different meanings!

We need some more specific terminology: with $(T,\cal T)$ denoting a topological space, where $T=\{1,...,n\}$, an $m-$\emph{system} is an open set $P$ of $\cal T$ such that $P-\{n\}$ has $m$ points. An \emph{upper} (resp., \emph{lower}) $m-$\emph{system} of $\cal T$ is an $m-$system which does not contain (resp., contains) $n$. A $k-$system $P$ of $\cal T$, with $k > 1$, is called an \emph{old} $k-$\emph{system} when each point $p$ of $P$ ($p\neq n$) is contained in some $m-$system where $m\leq k-1$, and is called a \emph{new} $k-$\emph{system} when at least one point $p$ of $P$ ($p\neq n$) is not contained in any $m-$system with $m\leq k-1$.

We denote by $\alpha ^{*}$ the smallest open set containing $\alpha \in T$ and we call $\alpha ^{*}$ the \emph{covering set} of $\alpha $; obviously, in a finite topological space, the covering sets completely determine the topology.
For a given value of $k$, we say that $\alpha (\neq n)$ is an \emph{old point} when  $\alpha ^{*}$ is an $m-$system for $m<k$; when $\alpha ^{*}$ is an $m-$system for $m=k$, we say that $\alpha $ is a \emph{new point}. In an \emph{old $k-$system}, all points distinct from $n$ are old. In a \emph{new $k-$system}, at least one point distinct from $n$ is new. Note that, for each $k$, there might be points which are neither old nor new, namely those points $\alpha $ such that $\alpha ^{*}$ is an $m-$system for $m>k$.

We write $m$-set to mean a set with $m$ points. It follows from the definitions that the existence of an open $m-$set in a quotient-space requires the existence in $(T,\cal T)$ of an open set with either $m$ points, all distinct from $n$ (as said above, we call it an \emph{upper $m-$system}), or $m+1$ points, one of them being $n$ (as said above, we call it a \emph{lower $m-$system}). More precisely: Let $A\subset T-\{n\}$ (hence $A \neq \{1, ..., n-1\}$); when $A$ is open in $(T, \cal T)$ but $A\cup \{n\}$ is not, the set $A$ is open in $Q^{i}$ if and only if $i\notin A$;
when $A\cup \{n\}$  is open in $(T, \cal T)$ but $A$ is not, the set $(A-\{ i\})\cup \{z\}$ is open in  $Q^{i}$ if and only if $i\in A$; when both $A$ and  $A\cup \{n\}$ are open in $(T,\cal T)$, $A$ is open in all spaces $Q^{i}$.
Usually, when referring to sets in  $Q^{i}$, where we should write $z$, we simply write $i$; no misunderstanding will arise. In fact, when $i,n\notin A$, we can naturally identify $A\subset T$ and $f(A)\subset Q^{i}$.

In a similar spirit, to adopt a shorter notation, when $A$, $B$, ... are subsets of the topological space $(T,\cal T)$, we write $AB$ to mean the union of the disjoint sets $A$ and $B$, and $An$ to mean the union of the disjoint sets $A$ and $\{n\}$. This shorter notation is very convenient to display the configurations introduced in the statements of Theorems 3, 4 and 5.
Moreover, sometimes we list the open sets of a topology by just writing their elements.

\section {Outline for a solution}

Our strategy to reconstruct a topological space $(T,\cal T)$ is to obtain its open sets or, more precisely, its $m-$systems, for increasing values of $m$. It is an iterative procedure: for each value of $k$, the $k$-iteration yields all $k$-systems.

\vspace{0.1cm}

We start with the $1-$systems of $\cal T$. We obtain them from the open $1-$sets (or singletons) of the quotient-spaces. The homeomorphism allows us to use the first integers $1, 2, ...$ as names for the points in the $1-$systems of $\cal T$ and then successive integers for the successive $m-$systems, as $m$ increases.

Suppose we have all $m-$systems for $m\leq k-1$. Let $X_{k-1}=\{1,...,x\}$ be the set of points distinct from $n$ and contained in them (obviously, for $k=1$ we have $X_{0}=\emptyset $).
Moreover, suppose that we can identify (we'll see how to do this in Section 7), for each value of $k$, the so-called \emph{old spaces} which are those obtained by identification of $n$ with each one of the points in $X_{k-1}$, the remaining spaces being called, for this value of $k$,  \emph{new spaces}.

For each value of $k$, the \emph{old} $k-$systems are those whose points have appeared in $k'-$systems with $k'<k$; the new ones may be
\emph{clean}, when they are covering sets of all its points (except possibly $n$ if they contain it) or \emph{mixed} when they contain \emph{old points}  (points in $m-$systems for $m\leq k-1$). For these $k-$systems,
we have to find
the respective \emph{configuration} which describes the type (upper, lower, paired) of the new $k-$systems of $\cal T$.
There are 13 possible configurations, and, for a certain value of $k$, we may obtain more than one acceptable configuration. We then resort to the \emph{clans}, families of open sets we introduce in Section 5, to try to identify the right configuration. At this point, after obtaining all $m-$systems for $m\leq k$, we form the set $X_{k}=X_{k-1}\cup Y_{k}$, where $Y_{k}=\{x+1, ..., w\}$ is the set of points whose covering set is a $k-$system, and we list the quotient-spaces $\ Q^{x+1}$, ..., $\ Q^{w}$ which become \emph{old spaces}.

\vspace{0.2cm}

Since $T$ is finite, the procedure terminates as soon as all covering sets (and, consequently, all open sets) of $T$ are formed, with just one exception, namely, the set $\{1,...,n-1\}$. In fact, if  $\{1,...,n-1\}$ is the union of open $k-$sets with $k<n-1$, then  $\{1,...,n-1\}$ is an open set of the reconstructed topology $\cal T$. Otherwise, it is optional to consider it open or not: note that the set of all elements of any quotient-space $Q$ is always the image, through the respective natural map $f: T \rightarrow Q$, of the open set $\{1,...,n\}= T$ of $\cal T$. Obviously, by the very definition of topology, nothing prevents us from taking also $\{1,...,n-1\}$ as open set of $\cal T$. See \emph{Examples 1 and 2}.

\section{The 13 configurations}

The following Lemmas give us some perception about relations which occur among the $k-$systems we deal with in this paper. Let again $A$, $B$, ... and $R$, $S$, ... be subsets of the topological space $(T,\cal T)$. Recall that the $k$-iteration  of the algorithm we are going to present is the procedure which allows us to obtain all open $k$-systems.

\begin{lema}
If, at the $k-$th iteration, $AR$ is a new upper $k-$system where $A$ contains only old points and $R$ contains only new points, then $A$ is open.
\end{lema}

\emph{Proof:} Let $A=\{\alpha _{1}, ..., \alpha _{i}\}$. For $j=1,...,i$, we have $\alpha _{j}^{*}\subseteq AR$. [Obviously, $n\notin \alpha _{j}^{*}$ otherwise $AR\cap \alpha _{j}^{*}$ would be open and would contain $\alpha _{j}$ but not $n$, a contradiction].
Moreover, since no point of $R$ is contained in an open set with less than $k$ points, we may write $\alpha _{j}^{*}\subseteq A$. Hence, $\cup _{j} \alpha _{j}^{*}=A$.

\begin{lema}
If, at the $k-$th iteration, $ARn$ and $BSn$ are new lower $k-$systems where $A$ and $B$ contain only old points whereas $R$ and $S$ contain only new points, then $An$ and $Bn$ are open.
\end{lema}

\emph{Proof:} As in Lemma 1, for each $\alpha _{j}$ in $A$ and each $\beta _{j}$ in $B$, we have $\alpha _{j}^{*}\subseteq An$ and $\beta _{j}^{*}\subseteq Bn$. Moreover, noting that $R$ and $S$ are disjoint, we get  $n^{*}\subseteq ARn\cap BSn=An\cap Bn$. Hence, $n^{*}\cup ( \cup _{j} \alpha _{j}^{*})=An$ and
$n^{*}\cup ( \cup _{j} \beta _{j}^{*})=Bn.$ Look that, whether or not $A=B$, the sets $R$ and $S$ if not equal must be disjoint, otherwise  take $x\in R\cap S$, hence $x\in AR\cap BS$, hence $x$ belongs to a $k'$-system with $k'<k$ which means that not all points of $R$ and $S$ are new, against the hypothesis. If $A\neq B$ and $R=S$, then for $x\in R$ we get $x\in AR\cap BR$, hence $x$ is an old point, against the hypothesis. When $A=B$ and $R=S$ we have the hypothesis of Lemma 3.

\begin{lema}
If, at the $k-$th iteration, $ARn$ is a new lower $k-$system where $A$ contains only old points and $R$ contains only new points, then at least one of the two sets $A$ or $An$ is open.
\end{lema}

\emph{Proof:} No $\alpha _{j}^{*}$ can contain any point of $R$ but it may contain $n$. Hence $\cup \alpha _{j}^{*}$ is $A$ or $An$.

\begin{lema}
If, at the $k-$iteration, $AR$ and $ARn$ as well as $BS$ and $BSn$ are new $k-$systems, with $R$ and $S$ distinct and containing only new points, and $A$ and $B$ containing only old points,  then $A$ and $An$ as well as $B$ and $Bn$ are also open sets.
\end{lema}

\emph{Proof:} By Lemma 1, $AR$ open implies $A$ open and $BS$ open implies $B$ open. As in Lemma 2, $R$ and $S$ being distinct they must be also disjoint. Moreover, $n^{*}\subseteq ARn\cap BSn=An\cap Bn$, hence $n^{*}\subseteq An$ and $n^{*}\subseteq Bn$ which implies $\cup \alpha _{j}^{*}\cup n^{*}=An$ and similarly, with an obvious notation, $\cup \beta _{j}^{*}\cup n^{*}=Bn$.

\begin{lema}
If $X$ and $Xn$ is a paired $k-$system and $Yn$ is a lower $k'-$system with $Y\subset X$, then $Y$ is open, that means, $Y$ and $Yn$ form a paired $k'-$system.
\end{lema}

\emph{Proof:} It is enough to notice that $Yn\cap X=Y$, hence $Y$ is open.

\begin{teo}
All $k-$systems whose points, distinct from $n$, are in $X_{k-1}$ (that is, are old points at the $k-$th iteration) may be obtained as unions of $k'-$systems ($k'<k$), except possibly $n^{*}$ (when it is also a $k-$system).
\end{teo}

\emph{Proof:} First, let $A=\{\alpha _{1}, ...,\alpha _{k}\}$ be an upper $k-$system and let $A\subseteq X_{k-1}$. This means that, for $i=1,...,k$, $\ \alpha _{i}^{*}$ is a $k'-$system ($k'<k$) and $n\notin \alpha _{i}^{*}$. Hence, $A=\cup _{i}\alpha _{i}^{*}$, which proves the assertion in this case.

Now take a lower $k-$system, say $An=\{\alpha _{1}, ...,\alpha _{k},n\}$ with $A\subseteq X_{k-1}$, as in the preceding case. Again we have $An=\alpha _{1}^{*}\cup ...\cup \alpha _{k}^{*}\cup n^{*}$ where each $\alpha _{i}^{*}$ is a $k'-$system ($k'<k$). Suppose that $n^{*}\neq An$. This means that $\{n\}\subseteq n^{*}\subset An$, that is to say, $n^{*}$ is also a $k'-$system with $k'<k$. This completes the proof.

\vspace{0.1cm}

The reasoning in this proof fails when $n^{*}=An$. In fact, we may have, for instance, $1^{*}=\{1\}$, $2^{*}=\{2\}$ as $1-$systems and $n^{*}=\{1,2,n\}$. Here $X_{1}=\{1,2\}$, $A=\{1,2\}\subseteq X_{1}$ and $n^{*}=An=\{1,2,n\}$ is a $2-$system but is not the union of $1-$systems.

\vspace{0.1cm}

The importance of this theorem lies in the fact that, in accordance with it, unless $n^{*}=An$ with $A\subseteq X_{k-1}$ and $|A|=k$, the number of new open $k-$sets in each quotient-space may be obtained when a list of all $k'-$systems, $k'<k$, is already known. Recall that an upper $k-$system $A=\{\alpha _{1},...,\alpha _{k}\}$ yields an open $k-$set in all quotient-spaces \emph{except} $Q^{\alpha _{1}}$, ..., $Q^{\alpha _{k}}$ and a lower $k-$system $An=\{\alpha _{1},...,\alpha _{k},n\}$ yields an open $k-$set in $Q^{\alpha _{1}}$, ..., $Q^{\alpha _{k}}$. The number of new open $k-$sets in the \emph{new quotient-spaces} may therefore be determined without taking into account whether or not $n^{*}=An$, with $A\subseteq X_{k-1}$; the number of new open $k-$sets in the \emph{old quotient-spaces} may be affected by only one unit.

\begin{teo}
For given $k$, the total number of new open $k-$sets in each new quotient-space differ by at most 2 units, that means, we can say it is $s$, $s+1$ or $s+2$, where $s\geq 0$ is an integer which depends on $k$.
\end{teo}

\emph{Proof:} We want to prove that the total number of new open $k$-sets in each new quotient-space differ by at most 2 units from one new space to another. First note that, if $A$ (with $n\notin A$) and $An$ are both $k-$systems (we then say they are paired), then, by the definitions, all quotient-spaces contain either $A$ or $(A-\{i\})\cup \{z\}$ as open $k-$set. Further, if $A$ is an upper, non-paired $k-$system, then $A$ appears as open $k-$set in all quotient-spaces except those which correspond to the points belonging to $A$ and, finally, if $An$ is a lower, non-paired $k-$system, then $(A-\{i\})\cup \{z\}$ appears as open $k-$set in the quotient-spaces $Q^{i}$ which correspond to the new points belonging to $A$.

The assertion now follows very easily:

Suppose that, with $s\geq 0$, there are $s+1$ (this means, at least one) new upper (non-paired) $k-$systems. Each new point appears either in only one $k-$system or in two $k-$systems which are paired, otherwise it would appear in a $k'-$system with $k'<k$. Consider the quotient-space $Q^{i}$. If $i$ is a new point which appears in an upper, non-paired $k-$system, then $Q^{i}$ contains $s$ new open $k-$sets. If $i$ appears in two $k-$systems which are paired or if $i$ does not appear in $k-$systems, then $Q^{i}$ contains $s+1$, new open $k-$sets; and if $i$ appears in a lower, non-paired $k-$system, then $Q^{i}$ contains $s+2$ new open $k-$sets.

Now suppose that there are no upper, (non-paired) $k$-systems. If $i$ appears in two $k-$systems which are paired or if $i$ does not appear in $k-$systems, then $Q^{i}$ contains, say, $p$ new open $k-$sets; and if $i$ appears in a lower, non-paired $k-$system, then $Q^{i}$ contains $p+1$ new open $k-$sets.
This completes the proof.

\vspace{0.1cm}

Since, for each $k$, the number of new open $k$-sets in the new quotient-spaces differ by at most 2 units from one new space to another, now, regardless of which quotient-spaces have $s$, $s+1$ or $s+2$ new open $k$-sets, we
denote by $\mu _{1}$, $\mu _{2}$ and $\mu _{3}$ the number of new quotient-spaces with  $s$, $s+1$ and $s+2$ new open $k-$sets, respectively, and we distinguish four cases:

\vspace{0.1cm}

Case 1:  $\mu _{1}\neq 0$, $\mu _{2}=\mu _{3}=0$;

Case 2:  $\mu _{1}\neq 0\neq \mu _{3}$, $\mu _{2}=0$;

Case 3:  $\mu _{1}\neq 0\neq \mu _{2}$, $\mu _{3}=0$;

Case 4:  $\mu _{1}\neq 0$, $\mu _{2}\neq 0$, $\mu _{3}\neq 0$.

\vspace{0.1cm}

Since we do not rule out $s=0$ in any case, we see that $s=0$ in Case 1 means that there are no new open $k-$sets in any new quotient-space.

Each one of these 4 cases may be yielded by one of several \emph{configurations}.
Theorems 3, 4 and 5 list configurations of Type 1, 2, 3 or 4 associated with Cases 1, 2, 3 or 4, respectively; for clarity, we shall use a graphical layout in which upper and lower systems are written in an upper and lower position, respectively.

\begin{teo}
In Case 1, the new $k-$systems form one of the following configurations, where, for any $i$, $A_{i}$ is a set of old points, $R_{i}$ is a set of new points, and $\mu_{1}$ is the total number of new points:

1a) $s+1$ upper $k-$systems, with each new point appearing in exactly one of the sets $R_{i}$, say

\hspace{1cm}$A_{1}R_{1}\ \ $; ... ;$\ A_{s+1}R_{s+1}\ \ $;

1b) $s$ pairs of $k-$systems (here, $\mu_{1}$ is the total number of new spaces, alias, all old and new spaces receive these $s$ open sets with $k$ points), say

\hspace{1cm}$A_{1}R_{1}\ \ \ $; ... ;$\ A_{s}R_{s}\ $;

\hspace{1cm}$A_{1}R_{1}n\ $; ... ;$\ A_{s}R_{s}n\ $;

1c) the empty configuration, that is, no new systems (here $s=0$);

1d) a set of lower $k-$systems, with each new point appearing in exactly one of the sets $R_{i}$, say

\hspace{1cm}$A_{1}R_{1}n\ \ $; ... ; $\ A_{x}R_{x}n\ \ $.

\end{teo}

Note that, in this theorem, 1c and 1d require $s=0$ and $s=1$, respectively; 1c may be considered as 1b for $s=0$; and in 1a, 1b and 1d, the sets $A_{i}$ of old points are arbitrary, that is, they may be empty or non-disjoint, for instance. Obviously, for $k=1$, all sets $A_{i}$ are empty.

\vspace{0.2cm}

\emph{Proof:} By an argument like the one used in the proof of Theorem 2, we see that, in the present hypothesis, there cannot exist a new upper $k-$system and a new lower $k-$system unless they are paired. Moreover, if two paired new $k-$systems exist, then all new $k-$systems must be paired. The possible configurations are therefore those indicated.

In 1a, $R_{1}\cup ... \cup R_{s+1}$ contains all new points, otherwise there would be quotient-spaces with $s+1$ new $k-$sets, which is a contradiction. Similarly, in 1d, $R_{1}\cup ... \cup R_{x}$ contains all new points, otherwise there would be quotient-spaces with no new $k-$set.  This completes the proof.

\vspace{0.2cm}

We deal with the configurations of Cases 2 and 4 together.

\begin{teo}
Using $A$, $B$ to denote sets of old points, and $R$, $S$ to denote sets of new points, the new $k-$systems in Cases 2 and 4 form one of the following configurations, the last one being only possible in Case 4:

2a) or 4a)  With $|R_{1}|+...+|R_{s+1}|=\mu _{1}$, $|R|=\mu _{3}$ and $\mu _{2}$ counting the number of remaining new quotient-spaces which is also the number of points which appear only in $k''-$systems with $k''>k$:

\hspace{1cm} $A_{1}R_{1}\ \ $; ... ; $\ A_{s+1}R_{s+1}\ \ $;

\hspace{5.2cm}$BRn\ \ $;

2b) or 4b) With $|R_{1}|+...+|R_{s+1}|=\mu _{1}$, $|S_{1}|+...+|S_{x}|=\mu _{3}$, $\mu _{2}$ again as in 2a) or 4a), and $B_{1}\cap ... \cap B_{x}\neq \emptyset $ (See that $B_{1}\cap ... \cap B_{x} \neq \emptyset $, otherwise $n^{*}=\{n\}$ and no upper non-paired system could exist):

\hspace{1cm} $A_{1}R_{1}\ \ $; ... ; $\ A_{s+1}R_{s+1}\ \ $;

\hspace{5.2cm} $B_{1}S_{1}n\ \ $; ... ;$\ B_{x}S_{x}n\ \ $;

4c) With $|R_{1}|+...+|R_{i-1}|=\mu _{1}$, $|R_{s+2}|+...+|R_{x}|=\mu _{3}$, and $\mu _{2}\geq |R_{i}|+...+|R_{s+1}|>0$ (Here, $A_{i}\cap ... \cap A_{x} \neq \emptyset $, same reason as in the preceding paragraph):

\hspace{0.7cm} $A_{1}R_{1}$;... ;$\ A_{i-1}R_{i-1}$;$\ A_{i}R_{i}\ $;...;$\ A_{s+1}R_{s+1}\ $

\hspace{4.2cm} $\ A_{i}R_{i}n$;...;$\ A_{s+1}R_{s+1}n$;$\ A_{s+2}R_{s+2}n$;...;$\ A_{x}R_{x}n$.

\end{teo}

\emph{Proof:} As a first remark, note that, for $k=1$, all sets $A$ and $B$ are empty, hence configurations 2b, 4b and 4c never occur. Otherwise,
as in the preceding theorem, we see that the new $k-$systems cannot be all upper, all lower or all paired. Moreover, we cannot have configurations with paired $k-$systems and lower, non-paired $k-$systems but without upper, non-paired $k-$systems. Similarly, we cannot have configurations with paired $k-$systems and upper, non-paired $k-$systems but without lower, non-paired $k-$systems. Hence, a configuration must have upper and lower $k-$systems. If it has also paired ones, then $\mu _{2}\neq 0$. If it has no paired ones, then it may be $\mu _{2}=0$ or $\mu _{2}\neq 0$;  in fact, in this hypothesis, $\mu _{2}\neq 0$ is the number of points which appear only in $k''-$systems with $k''>k$.
 This completes the proof.

\begin{teo}
Using again the same notation as in the preceding theorems, the new $k-$systems in Case 3 form one of the following configurations:

3a)  With $|R_{1}|+...+|R_{s+1}|=\mu _{1}$  and $\mu _{2}> 0$ counting the number of remaining new quotient-spaces, which is also the number of points which appear only in $k''-$systems with $k''>k$:

\hspace{1cm} $A_{1}R_{1}\ $; ... ;$\ A_{s+1}R_{s+1}\ \ $;

3b) With $|R_{1}|+...+|R_{i-1}|=\mu _{1}$, $\mu _{2}\geq |R_{i}|+...+|R_{s+1}|\ $, $A_{i}\cap ... \cap A_{s+1} \neq \emptyset $:

\hspace{0.9cm} $A_{1}R_{1}$;... ;$\ A_{i-1}R_{i-1}$;$\ A_{i}R_{i}\ $;...;$\ A_{s+1}R_{s+1}\ $

\hspace{4.4cm} $\ A_{i}R_{i}n$;...;$\ A_{s+1}R_{s+1}n$\ ;

3c) With $|R_{1}|+...+|R_{x}|=\mu _{2}$ and $\mu_{1}$ counting the number of the remaining quotient-spaces among the new ones:

\hspace{4,9cm} $A_{x+1}R_{x+1}\ $; ... ;$\ A_{x+s}R_{x+s}\ \ $;

\hspace{1cm} $A_{1}R_{1}n$; ... ; $\ A_{x}R_{x}n$; $\ A_{x+1}R_{x+1}n\ $; ... ;$\ A_{x+s}R_{x+s}n$;

3d) With $|R_{1}|+...+|R_{x}|=\mu _{2}$  and $\mu _{1}> 0$ counting the number of points which appear only in $k''-$systems with $k''>k$ which is the same as the number of remaining quotient-spaces among the new ones:

\hspace{1cm} $A_{1}R_{1}n$; ... ;$\ A_{x}R_{x}n\ $.

\end{teo}

Note that, in the statement of this theorem, 3d requires $s=0$; 3a, for $\mu_{2}=0$, is 1a; 3d, for $\mu _{1}=0$, is 1d. Note also that, for $k=1$, we cannot have, in 3b, more than one paired system.

\emph{Proof:} In this case, we cannot have one lower, non-paired $k-$system together with one upper, non-paired $k-$system in the same configuration; in fact, if $r_{1}$ is a new point which appears in a lower, non-paired $k-$system, and $r_{2}$ is a new point which appears in an upper, non-paired $k-$system, then $Q^{r_{1}}$ has two more open $k$-sets than $Q^{r_{2}}$, which is impossible in Case 3. Hence, the possible configurations are those indicated and the theorem is proved.

\vspace{0.2cm}

Let us illustrate what we have just said with an example.

\vspace{0.1cm}

\emph{Example 3:} suppose we are given three quotient-spaces: $Q_{1}= \{a,b,c\}$ with $\emptyset $, $\{a\}$, $\{a,b,c\}$ as open sets, $Q_{2}= \{x,y,z\}$ with $\emptyset $, $\{x\}$, $\{x,y,z\}$ as open sets, $Q_{3}= \{g,h,j\}$ with $\emptyset $, $\{g\}$, $\{h\}$, $\{g,h\}$, $\{g,h,j\}$ as open sets.

For $k=1$, we have $\mu_{1}=2$ spaces with $s=1$ singleton each and $\mu_{2}=1$ space with $s+1=2$ singletons. It is \emph{configuration 3a} with $\{1\}$ and $\{2\}$ as upper $1-$systems. As old spaces, choose [we'll see in Section 7 how to do these choices] $Q^{1}: \emptyset,\ \{a\};\ \{a,b,c\}$ and $Q^{2}: \emptyset ,\ \{x\},\ \{x,y,z\}$. The space $Q_{3}$, with $\emptyset $, $\{g\}$, $\{h\}$, $\{g,h\}$, $\{g,h,j\}$ as open sets, remains as only new space.

For $k=2$, we have one old open $2-$set $\{g,h\}=\{1,2\}$; no new open $2-$set in the new space. This is  \emph{configuration 1c}.

For $k=3$, we have in the same space $\mu_{1}=1$ and $s=1$ with $\{g,h,j\}=\{1,2,3\}$ as a unique new open $3-$set. It is \emph{configuration 1b} because $\{1,2,3,4\}=T$ is also open.

Our conclusion: these were the quotient spaces of $(T, \cal T)$ where $T=\{1,2,3,4\}$, $n=4$ and
$\emptyset $, $\{1\}$, $\{2\}$, $\{1,2\}$, $\{1,2,3\}$, $\{1,2,3,4\}$ the open sets of $\cal T$.

\section{Clans: definitions and their role}

Clans and configurations are the basic tools for the reconstruction of $\cal T$, or, more precisely, for finding solutions to our \emph{Problem}.
Let us define \emph{clan of open sets} and \emph{clan of covering sets}.

We resort to Graph Theory (see \cite{Harary} or \cite{Livromeu})) and first we associate, to a finite topology, a digraph $G$ whose vertices are the open sets of the topology and whose arcs are defined as follows: there is an arc from $\beta $ to $\alpha $, written  $(\beta , \alpha )$, when $\beta $ is properly contained in $\alpha $ and there is no $\gamma $ such that $\beta \subset \gamma \subset \alpha $. There is an exception to this rule: $A$ will not be connected to $An$.

A \emph{clan (of open sets)} is a \emph{connected sub-digraph} of $G$; we don't say a \emph{connected component} of $G$  because we don't require maximality and we exclude sub-digraphs with the vertex associated with $\emptyset $.
The designation \emph{clan} evokes our concept of a family.
When there is a directed path from $\beta $ to $\alpha $, we say that $\alpha $ is a descendant of $\beta $ and $\beta $ an ancestor of $\alpha $. When there is one single arc $(\beta , \alpha )$, we say that $\alpha $ is in the generation following the generation of $\beta $.
A vertex with no ancestor is called a \emph{root}. Obviously, $\emptyset $ would be an ancestor of all vertices of the sub-digraph, but we have excluded it in the definition of a clan. As regards the whole set where the topology is defined, it is obvious that it is a descendant of every open set of the topology; when it is the only descendant of the root (or it is the root itself!), the clan is called \emph{trivial}.

Two clans are said to be \emph{isomorphic} when their associated digraphs are isomorphic and corresponding vertices in both digraphs are open sets with the same number of points. When only the first of these two conditions is met the clans are said to be \emph{similar}.

For a \emph{clan of covering sets} we can say the same we just said for clans of open sets, just don't worry about $\emptyset $ because no covering set is empty.

Note that a clan may have more than one root, and a set may belong to more than one clan, each one of them with a distinct root. For instance, with $T=\{1,...,8\}$, the sets $1^{*}=\{1\}$, $2^{*}=\{2\}$, $3^{*}=4^{*}=\{1,2,3,4\}$, $5^{*}=\{1,5\}$, $6^{*}=7^{*}=\{2,6,7\}$ and $8^{*}=T$ form a clan. We can as well recognize here two clans which are distinct but not disjoint: one is rooted at $\{1\}$ with $\{1,5\}$, $\{1,2,3,4\}$ and $T$ as sets, another is rooted at $\{2\}$ with $\{2,6,7\}$, $\{1,2,3,4\}$ and $T$ as sets. They are \emph{similar} but not \emph{isomorphic} because the set $\{1,5\}$ has not the same number of points as the set $\{2,6,7\}$.

Given a clan, a \emph{subclan} is formed by one of its sets in the role of root and the respective descendants.

\vspace{0.1cm}

The following observations are immediate consequences of the definitions and the Lemmas given in Section 4:

\vspace{0.1cm}

\emph{Observation 1}: A clan with just one root in a lower non-paired system contains only lower non-paired systems.

\emph{Observation 2}: A clan sometimes said to have two roots in a paired system may be better understood as two clans with two distinct roots.

\emph{Observation 3}: A clan with just one root in an upper non-paired system may contain all kinds of systems.

\emph{Observation 4}: A clan rooted on an upper non-paired system and with no lower non-paired system appears entirely in $Q^{j}$ if and only if $j$ is not a point of any system of the clan.

\emph{Observation 5}: A clan rooted on a lower non-paired system yields a \emph{similar} clan in $Q^{j}$ if and only if $j$ is one point of the root.

\vspace{0.1cm}

\textbf{A list of clans:} As a very useful auxiliary step to rebuild the topology $\cal T$ defined on $T$, look at the given quotient-spaces and identify as many clans as possible in the topology of each  quotient-space. Note that the clans of $\cal T$ on $T$ yield clans in the quotient-spaces.

Do not forget that a point represented by a certain symbol in one space may be represented by a different symbol in another space. For instance, a clan like $5; 56; 567$ of $\cal T$ may yield in $Q^{1}$,  $\ p; pq; pqr$ and, in $Q^{2}$, $\ d; df; dfg$. What really matters is how the sets of the clan relate among them.

The names of the symbols may, and usually do vary, from space to space. But we'll be able to distinguish clans of lower non-paired systems, clans of upper non-paired systems and pairs of clans with their roots and all sets in paired systems. Keep all these clans at hand while looking at the configurations.

It is also important to notice that some apparently isomorphic clans may be of different types, that is, one lower, the other upper. Two disjoint upper isomorphic clans appear together in spaces $Q^{j}$ where $j$ does not belong to their sets; two disjoint lower isomorphic clans never appear together; when one is upper and the other is lower, they may appear together in $Q^{j}$ where $j$ is any point of the root of the lower one: as pointed out above, the lower clan yields a \emph{similar} one in $Q^{j}$.
Clans with their roots in a paired system yield a similar clan in all spaces.

Keep in mind that in a clan of lower systems, the systems with more points appear in more spaces; and in a clan of upper systems, the systems with more points appear in fewer spaces.

\section{The algorithm: role of the configurations}

In this and the next section, we describe the steps of our algorithm. We keep in mind Theorems 2, 3, 4 and 5 of Section 4 and if, at any step, for any value of $k$, the conditions they presume are not respected, then no solution exists; the algorithm, if programmed for a computer, shall give us a message that no solution exists.

Let us see what happens when the above mentioned conditions are always respected.
As suggested in Section 3, we take the given quotient-spaces and, for each $k$, using the $k-$sets they contain, we try to find the configurations which might have given rise to such sets.

\vspace{0.2cm}

\textbf{\emph{To begin with, consider the case where $k=1$}}. This is the base case of our iteration procedure. Since we look for $\cal T$ to within homeomorphism, let us choose as names for the points in the $1-$systems the first natural numbers $1, 2, ..., $; by other words, give to the points of T which constitute covering sets with just one point or one point and $n$, the first integer numbers.

\vspace{0.1cm}

Suppose all spaces have the same number $s$ of $1-$sets or singletons. These are configurations of type 1. If $s=0$, we have \emph{configuration 1c}; if $s\neq 1\wedge s\neq n-2$, we have \emph{configuration 1b}; if $s=1$ we have \emph{configuration 1b or 1d} and we can make the right choice by checking whether there are new open $2-$sets (or $k'-$sets for $k'\geq 2$), which means \emph{1b}, or there are none, which means \emph{1d}; if $s=n-2$ we have \emph{configuration 1b or 1a} and here again, when there are new open $2-$sets (or $k'-$sets for $k'\geq 2$), we choose \emph{1b}; if there are no new open sets (or $k'-$sets for $k'\geq 2$), we may have two configurations which apparently yield two solutions, namely, $\emptyset $, $\{1\}$, $\{2\}$, ..., $\{n-1\}$, $n^{*}=T$ and $\emptyset $, $\{1\}$, $\{1,n\}$, $\{2\}$, $\{2,n\}$,..., $\{n-2\}$, $\{n-2,n\}$, $(n-1)^{*}=T$, $n^{*}=\{n\}$, but note that these two topologies are homeomorphic.

\vspace{0.1cm}

As regards configurations of type 2 or 4, the only possibilities for $k=1$ are as follows: we may have $\mu_{1}$ spaces with $s$ singletons and $\mu_{3}=1$ space with $s+2$ singletons, that is, \emph{configuration 2a} when $s+2=n-1$; or \emph{configuration 4a} when $s+2< n-1$, a case where there will be $\mu_{2}=n-1-(s+2)$ spaces with $s+1$ singletons. In both cases we take $\{1\}$, $\{2\}$, ... $\{s+1\}$, $\{s+2,n\}$, as open $1-$systems.

\vspace{0.1cm}

Finally, look at configurations of type 3. \emph{Configuration 3a} yields $\mu_{1}=s+1$ spaces with $s$ singletons and the remaining $\mu_{2}=n-1-(s+1)$ spaces with $s+1$ singletons. \emph{Configuration 3d} yields $\mu_{2}=x$ spaces with 1 singleton and the remaining $\mu_{1}=n-1-x$ spaces with no singleton. \emph{Configuration 3b}, when $k=1$, admits only one paired $1-$system, hence yields $\mu_{1}=s$ spaces with $s$ singletons and $\mu_{2}=n-1-s$ spaces with $s+1$ singletons. \emph{Configuration 3c} yields $\mu_{2}=x$ spaces with $s+1$ singletons and the remaining $\mu_{1}=n-1-x$ spaces with $s$ singletons.

\vspace{0.1cm}

Look at the following synoptic table  where we give, for each configuration of Type 3, the number of spaces and the number of singletons in each space:

\vspace{0.2cm}

\hspace{0.0cm}Config.\hspace{0.5cm} Spaces\hspace{0.8cm}    Singletons\hspace{1.3cm}  Spaces\hspace{1.5cm}     Singletons

\hspace{+0.3cm}$3a$ \hspace{0.5cm} $\ \mu_{1}=s+1$ \hspace{1cm}$s$ \hspace{0.8cm} $\mu_{2}=n-1-(s+1)$   \hspace{0.8cm} $s+1$

\hspace{+0.3cm}$3d$ \hspace{0.2cm} $\ \mu_{1}=n-1-x$ \hspace{0.6cm}$0$ \hspace{1.9cm} $\mu_{2}=x$   \hspace{2.2cm} $1$

\hspace{+0.3cm}$3b$ \hspace{0.7cm} $\ \mu_{1}=s$ \hspace{1.6cm}$s$ \hspace{1cm} $\mu_{2}=n-1-s$   \hspace{1.5cm} $s+1$

\hspace{+0.3cm}$3c$ \hspace{0.3cm} $\ \mu_{1}=n-1-x$ \hspace{0.6cm}$s$ \hspace{1.9cm} $\mu_{2}=x$   \hspace{1.9cm} $s+1$

\vspace{0.1cm}

Note that this table holds only for $k=1$. For $k>1$, replace $n-1$ by $n-1-t_{k}$ with $t_{k}$ counting the number of old spaces at the $k$ iteration (obviously, $t_{1}=0$), and the word "singletons" by "new open $k-$sets".

See now how we can choose the right configuration:

3a versus 3d: To have $s=0$ in 3a, we need $\mu_{1}=1$ and $\mu_{2}=n-2$, that is, one space with no singleton and the remaining spaces with 1 singleton; in 3d, there is only one space with 1 singleton. Hence an ambiguity appears only when $n=3$. In fact, when $T=\{1,2,3\}$ and $n=3$, the pair of quotient-spaces $Q': \emptyset , \{x\}, \{x,y\}$ and $Q'': \emptyset , \{u,v\}$ may be yielded by the non-homeomorphic topologies ${\cal T}_{1}: \emptyset , \{1\}, T$ (configuration 3a) or ${\cal T}_{2}: \emptyset , \{1,3\}, T$ (configuration 3d). See that the clans of the quotient-spaces are trivial.

3a versus 3b: the number of spaces with $s$ singletons is $s+1$ in 3a and is $s$ in 3b.

3a versus 3c: in 3a, all spaces with $s$ singletons have $C^{s}_{2}$ doublets (or $2-$sets) formed by their singletons ($C^{s}_{2}$ denotes the binary combinatorial coefficient, and for $s=1$ we have $C^{1}_{2}=0$); in 3c, $s$ of the $\mu_{1}$ spaces with $s$ singletons have $C^{s}_{2}$ doublets formed by their singletons plus $x$ doublets formed by one of their $s$ singletons with each one of $x$ symbols.

3d versus 3b: in 3b, $\mu_{1}=0$ spaces with 0 singletons would lead to a Type 1 configuration.

3d versus 3c: configuration 3c becomes 3d if some spaces have 0 singletons.

3b versus 3c: keep in mind that, for $k=1$ configuration 3b can have only one paired system; now, if $\mu_{1}>s$ then choose 3c; if  $\mu_{1}=s$, then check the spaces with $s$ singletons: like what we saw in 3a versus 3c, in 3b the spaces with $s$ singletons have $C^{s}_{2}$ doublets (or $2-$sets) formed by their singletons; in 3c, $s$ of the $\mu_{1}$ spaces with $s$ singletons have $C^{s}_{2}$ doublets formed by their singletons plus $x$ doublets formed by one of their $s$ singletons with each one of $x$ symbols.

\vspace{0.2cm}

\textbf{\emph{Now consider the cases where $k>1$}.} Assign successive integers to the new points in the $k-$systems for $k=2,3,...$. When no $1-$systems exist, start here with $2-$systems or, more generally, with the smallest systems.

\vspace{0.2cm}

\textbf{Configurations of type 1.} When all new spaces have the same number, say $s$, of new $k-$sets, the possible configurations are 1a, 1b, 1c or 1d. Here we look separately at cases \textbf{A} ($s=0$), \textbf{B} ($s=1$) and \textbf{C} ($s>1$).

\vspace{0.2cm}

\textbf{A:} When $s=0$, two configurations are possible: 1a (with $s+1=1$ and hence $R_{1}$ containing all new points) or 1c. If, in 1a, we have only $A_{1}R_{1}$, that is, $s+1=1$, hence $s=0$, then we get $\mu_{1}=|R_{1}|$ new spaces with 0 new open $k$-sets, which means that the new $k-$system does not yield any open $k$-set in any new space. How can we distinguish now 1a from 1c? If $A_{1}R_{1}\neq \{1, ..., n-1\}$, then for $j\notin A_{1}R_{1}$, the space $Q^{j}$ contains $A_{1}R_{1}$. Of course $Q^{j}$ is an old space, otherwise the configuration was not 1a, because not all the new spaces would have the same number of new $k-$systems.
If $A_{1}R_{1}=\{1, ..., n-1\}$, then there is no need for further discussion, because all $q$-spaces have $\{1, ..., n-1\}$ as open set, either new, as the image $f(T)$, or old, an old $k$-set which is the union of open $k'$-sets for values of $k'< k$. By other words, for $k=n-1$, configuration 1c never occurs. See \emph{Example 7}, Section 9, for $k=2$ and $k=n-1=3$.

\vspace{0.2cm}

\textbf{B:} When $s=1$, we may have three configurations: 1a, 1b or 1d.
In 1a, we have two upper $k-$systems $A_{1}R_{1}$ and $A_{2}R_{2}$; in 1b, the paired $k-$systems $A_{1}R_{1}$ and $A_{1}R_{1}n$; in 1d, $x$ lower $k-$systems $A_{1}R_{1}n$, ..., $A_{x}R_{x}n$, with $1\leq x\leq n-1$.

This is a situation where we may have several non-homeomorphic topologies with the same quotient-spaces: as an example, take $T=\{1,2,3,4,5\}$, $n=5$ and four quotient-spaces $Q_{1}=Q_{2}=Q_{3}=Q_{4}: \emptyset, \{u,v\}, \{u,v,x,y\}$. For $k=1$ we have configuration 1c, but for $k=2$ we may choose configuration 1a, which yields ${\cal T}_{1}: \emptyset , \{1,2\}, \{3,4\}, \{1,2,3,4\}, T$ with $n^{*}=T$, configuration 1b, which yields ${\cal T}_{2}: \emptyset , \{1,2\}, \{1,2,n\}, T$ with $n^{*}=\{1,2,n\}$, or configuration 1d, which yields ${\cal T}_{3}: \emptyset , \{1,2,n\}, \{3,4,n\}, \{n\}, T$ with $n^{*}=\{n\}$. This is another example where the clans of the quotient-spaces are all trivial.

\vspace{0.1cm}

Sometimes we may be forced to make a choice. Look that,
in 1a, $|R_{1}|+|R_{2}|=\mu _{1}$ (all new points are in $R_{1}\cup R_{2}$, hence $|R_{1}|, |R_{2}|<\mu_{1}$) and $|A_{1}R_{1}|=|A_{2}R_{2}|=k$. When not all $k-$sets in the new spaces have the same number of new points, that means, when $|R_{1}\neq |R_{2}|$, then 1b is excluded. We are left with 1a or 1d. Choose 1a when $|R_{1}|$ new spaces contain $A_{2}R_{2}$ and $|R_{2}|$ new spaces contain $A_{1}R_{1}$; choose 1d when $|R_{2}|$ new spaces contain $A_{2}R_{2}$ and $|R_{1}|$ new spaces contain $A_{1}R_{1}$. Another situation when we may exclude 1b is when there is an old point $j$ such that $j\notin A_{1}R_{1}\cup A_{2}R_{2}$; then there is an old space $Q^{j}$ which contains both $A_{1}R_{1}$ and $A_{2}R_{2}$.

\vspace{0.2cm}

\textbf{C:} When $s>1$, we may have two configurations: 1a or 1b.

In 1b, all spaces (old and new) contain $s$ new open $k-$sets, but, obviously, not necessarily all new points.  In 1a, $R_{1}\cup ... \cup R_{s+1}$ contains all new points. In fact, if $x$ is a new point and $x\notin R_{1}\cup ... \cup R_{s+1}$, then $Q^{x}$ contains the $s+1$ new open $k-$sets, a contradiction. If there is some $x\notin A_{1}R_{1}\cup ... \cup A_{s+1}R_{s+1}$, then $Q^{x}$ has $s+1$ new open $k-$sets and it contains $A_{1}R_{1}$,..., $A_{s+1}R_{s+1}$. Of course $Q^{x}$ must be an old space; indeed, if $Q^{x}$ were a new space, then the configuration was not 1a, because not all the new spaces would have the same number of new open $k-$sets.

If no point $x$ exists such that $x\notin A_{1}R_{1}\cup ... \cup A_{s+1}R_{s+1}$, then check whether $|R_{1}|=...=|R_{s+1}|$. If these equalities don't hold, then choose 1a. If they hold, we may choose 1a or 1b. For 1a, $\forall i: |R_{i}|=r/(s+1)$. For 1b, when all new points are in new $k-$sets, then $\forall i: |R_{i}|=r/s$; when not all new points are in the new $k-$sets, we may have $\forall i: |R_{i}|=r/(s+1)$, but in this case, let $x$ be one point that is not in the new $k-$sets: to allow us the choice of 1b, $x$ must show up in the quotient-spaces through its covering $k'-$set, $k'>k$ unless $x^{*}=T-\{n\}$. This is a case where we may have 2 solutions (See \emph{Example 8} in Section 9).

\vspace{0.2cm}

\textbf{Configurations of type 2 or type 4.}
Let us count the number of new open $k-$sets in each new quotient-space. We distinguish configurations 2 from configurations 4 by the simple fact that, in configurations 2, some new spaces have $s$, others have $s+2$, but no one has $s+1$ new open $k-$sets, and in configurations 4, besides those with $s$ and $s+2$, there are also new spaces with $s+1$ new open $k-$sets.

\vspace{0.1cm}

The way to distinguish between configurations 2a and 2b is an immediate consequence of the definitions: in 2a, the spaces with $s+2$ sets exhibit all the new points in their new sets; in 2b, no space with $s+2$ sets exhibits all the new points in its new sets, the reason being that in each one of these spaces only one of the two or more lower systems of 2b will be present. Don't forget that the new points in each new space are those which appear only in open sets with $k$ or more points and, of course, never in sets with $k'<k$ points.

As regards configurations 4, it is important to register the points which have appeared as elements of open sets with $k''<k$ elements and count them; count as well those which now appear in new spaces as elements of open sets with $k$ points; this allows us to know how many are the remaining points, that is, those which appear in open sets with $k'>k$ elements. After doing this, the distinction between 4a and 4b is similar to what we did for 2a and 2b.

\vspace{0.2cm}

\textbf{Configurations of type 3.} It remains to consider the case when $\mu_{1}$ new spaces have $s$ new open $k-$sets and $\mu_{2}$ have $s+1$ new open $k-$sets, that is to say, configurations 3. We look separately at cases \textbf{A} ($s=0$) and \textbf{B} ($s\geq 1$).

\vspace{0.1cm}

\textbf{A:} When $s=0$, configurations 3b and 3c cannot occur; in these configurations, all new spaces have new systems. We may have 3a with $A_{1}R_{1}$ as only $k-$system, hence $\mu_{1}=|R_{1}|$ spaces with no new open $k-$set and all other new spaces with one new open $k-$set; or we may have 3d with $\mu_{2}=|R_{1}|+ ... +|R_{x}|$ new spaces with one new open $k-$set and the remaining $\mu_{1}$ new spaces with no new open $k-$set.

To distinguish 3a from 3d (when $s=0$), we may usually resort to the number of new points in the new spaces with one new $k-$system. If this number is not the same for all these spaces, choose 3d. If it is the same, more has to be done to distinguish 3d from 3a: we then resort to the values of $\mu_{1}$ and $\mu_{2}$. However, when  $\mu_{1}=\mu_{2}$, it may be impossible to make the distinction, even when we try to resort to the clans.

\vspace{0.1cm}

Recall the topologies ${\cal T}_{1}$ and ${\cal T}_{2}$ of \emph{Example 1}. All their clans are trivial. Non-trivial clans can lead to a choice between 3a and 3d, as seen below:

\vspace{0.1cm}

\emph{Example 4:} Look at the following topologies on the set $T=\{1,2,3,4,5\}$: ${\cal T}_{3}$ with $\emptyset $, T, $\{1,2,5\}$ and  $\{1,2,3,5\}$ as open sets; ${\cal T}_{4}$ with $\emptyset $, T, $\{1,2\}$ and  $\{1,2,3\}$ as open sets; ${\cal T}_{5}$ with $\emptyset $, T, $\{1,2\}$, $\{1,2,3\}$  and  $\{1,2,4\}$ as open sets.

\vspace{0.1cm}

In these topologies, one of the clans is rooted at a lower $2-$system and it forces us to choose configuration 3d; two other clans are rooted at an upper $2-$system and they force us to choose configuration 3a.

This can be recognized when we look at the spaces as they are given to us.

For $k=1$, the configuration is 1c for all topologies.

For $k=2$, all these topologies have $\mu_{1}=2$ spaces with $s=0$ new open $2-$sets and $\mu_{2}=2$ spaces with $s=1$ new open $2-$set. To decide between 3a (with $12$ as only $2-$system) and 3d (with $12n= 125$ as only $2-$system), we have to look at the clans. In ${\cal T}_{3}$ we recognize a clan whose sets with more points appear in more spaces; in ${\cal T}_{4}$ we recognize a clan whose sets with more points appear in fewer spaces. In ${\cal T}_{5}$ the spaces with $\alpha \beta; \alpha \beta \gamma $ as open sets are not produced by a lower clan; in some space the three systems of such a clan would yield open $3-$sets; hence we have a clan with three upper systems.
By other words, in ${\cal T}_{5}$ we have two spaces both with $\alpha \beta ; \alpha \beta \gamma $ as open sets meaning $\{1,2\},\ \{1,2,3\},\ \{1,2,4\}$ as open sets in the reconstructed topology; in ${\cal T}_{3}$ we have two spaces with $\alpha \beta; \alpha \beta \gamma $ plus one space with $\lambda \mu \nu$ as open sets meaning $\{1,2,5\},\ \{1,2,3,5\}$ as open sets in the reconstructed topology; in ${\cal T}_{4}$ we have one space with $\alpha \beta ; \alpha \beta \gamma $ plus one space with $\varphi \psi$ as open sets meaning $\{1,2\},\ \{1,2,3\}$ as open sets in the reconstructed topology. After having the $2-$systems, the old spaces $Q^{1}$ and $Q^{2}$  will be: for ${\cal T}_{1}$ (of \emph{Example 1}), ${\cal T}_{4}$ and ${\cal T}_{5}$ the spaces with no new open $2-$set; for ${\cal T}_{2}$ (of \emph{Example 1}) and ${\cal T}_{3}$, the spaces with one new open $2-$set.

For $k=3$, we need no further observation, but it may be curious to verify that in ${\cal T}_{1}$ and ${\cal T}_{2}$ (both of \emph{Example 1}) no new open $3-$set shows up in the new spaces, hence we have configuration 1c; in ${\cal T}_{3}$ and ${\cal T}_{4}$, we have $\mu _{1}=1$ new space with $s=0$ new open $3-$set and $\mu_{2}=1$ new space with $s+1=1$ new open $3-$set, which means again that, if we do not look at the clans, we would be able to choose 3a or 3d; finally, in ${\cal T}_{5}$, we have $\mu_{1}=2$ spaces with $s=1$ new open $3-$set, hence the configuration to be chosen is 1a.

For $k=4$, ${\cal T}_{1}$ and ${\cal T}_{2}$ (both of \emph{Example 1}) have two new spaces and $\alpha \beta \gamma \delta$ is the new open $4-$set which does appear in every space: it corresponds to configuration 1b; in ${\cal T}_{3}$ and ${\cal T}_{4}$ we have one new space and configuration 1b, again; finally, in ${\cal T}_{5}$ there is no new space for $k=4$ and $\{1,2,3,4\}$ is an old $4-$system.

\vspace{0.2cm}

\textbf{B:} When $s\geq 1$, configuration 3d cannot occur, but 3a, 3b and 3c are possible. Let us present a few features which help us choose the right configuration. As we pointed out, in a synoptic table for Type 3 and $k>1$, instead of $n-1$ as in the table for $k=1$ we should write now $n-1-t_{k}$ with $t_{k}$ counting the old spaces for $k$.

Suppose we know which spaces are new and that we recognize the new open $k-$sets. Let $p$ be the number of new points in a space with $s+1$ new open $k-$sets. If not all new spaces with $s+1$ new open sets have the same number $p$ of new points, then we are done: we choose configuration 3c.

Suppose now that all new spaces with $s+1$ new open sets have the same number $p$ of new points. If $\mu_{1}\neq p$, configuration 3a is excluded; we can have only 3b (where $\mu_{1}<p$) or 3c. If $\mu_{1}=p$ configuration 3b is excluded; we can have only 3a or 3c.

Let $p'$ be the number of new points in a space with $s$ new open $k-$sets. If $p'$ is not the same for the $\mu_{1}$ spaces with $s$ new sets, then the configuration cannot be 3c, it must be 3a or 3b. Hence we can say: when $p'$ is not the same for all spaces with $s$ systems, we choose configuration 3a if $\mu_{1}=p$, and configuration 3b if $\mu_{1}<p$.

At this point we see that we have to find ways to distinguish 3a from 3c (when $\mu_{1}=p$  for all new spaces with $s+1$ new open $k-$sets and $p'$ is the same for all new spaces with $s$ new open $k-$sets); or 3b from 3c (when $\mu_{1}<p$, $p$ is the same for all new spaces with $s+1$ new open $k-$sets and $p'$ is the same for all new spaces with $s$ new open $k-$sets).

To distinguish 3a from 3c, look at the new spaces with $s$ new open $k-$sets: for configuration 3a, we have $C^{s}_{2}$ combinations of these sets and they appear together with the $s$ sets, in such spaces; for configuration 3c, there will be, among the $\mu_{1}$ new spaces with $s$ new $k-$sets, $s$ spaces with $C^{s}_{2}$ combinations of the $s$ sets plus combinations of each one of the sets in lower non-paired systems with one of the sets in the paired systems.

In fact, in this latter case the union of one $k-$set associated with a lower, non-paired system, say $A_{1}R_{1}n$ with the $k-$set associated with the lower set of a paired system, say $A_{x+1}R_{x+1}n$, yields a set which appears in $Q^{\beta }$ for $\beta \in R_{x+1}$ and for $\beta \in R_{1}$.

\vspace{0.2cm}

To distinguish 3b from 3c, check the unions of the new open $k-$sets associated to paired $k-$systems with the new open $k-$sets associated with upper $k-$systems in configuration 3b, or with lower $k-$systems in configuration 3c. In 3b, these unions yield paired systems: for instance, $A_{1}R_{1}\cup A_{i}R_{i}$ and $A_{1}R_{1}\cup A_{i}R_{i}n$ yield a paired system, hence it appears in all spaces. In 3c,
$A_{1}R_{1}n\cup A_{x+1}R_{x+1}=A_{1}R_{1}n\cup A_{x+1}R_{x+1}n$, a lower system which does not appear in all spaces, unless
$A_{1}R_{1}\cup A_{x+1}R_{x+1}=\{1, ..., n-1\}$, as seen below:

\emph{Example 5:} Take $T=\{1,2,3\}$ with $n=3$ and quotient-spaces $Q_{1}$ and $Q_{2}$ with topologies $U_{1}$ whose open sets are $\emptyset $, $a, b, ab$ and $U_{2}$ whose open sets are $\emptyset $, $u, uv$. If seen as configuration 3b (for $k=1$) we reconstruct $\cal T'$ with $\emptyset $, T, $\{1\}$, $\{2\}$, $\{2,3\}$ as open sets; if seen as configuration 3c (for $k=1$) we reconstruct $\cal T''$ with $\emptyset $, T, $\{1,3\}$, $\{2\}$, $\{2,3\}$ as open sets. Nonetheless, it is easy to see that
$\cal T'$ and $\cal T''$ are homeomorphic: associate to points $1, 2, 3$ of $\cal T'$, points $2, 3, 1$, respectively, of $\cal T''$.

In the cases where we have found more than one solution, the quotient-spaces have only trivial clans (or, more precisely, trivial subclans) with $k-$sets as roots.

\section{The algorithm: choosing the old spaces}

As we already said, for each value of $k$ we have to identify the \emph{old spaces}, before moving ahead. However we believe it is better to explain, in this separate section, how to do this identification. Recall that, for each value of $k$, some $k-$systems are entirely formed by old points, that is, points which appear in $m-$systems for $m\leq k-1$; a new $k-$system contains at least one new point, that means a point which has not appeared in any $m-$system for $m\leq k-1$. Such points are new for $k$, but become old, and the spaces obtained by identification of $n$ with each one of them will become \emph{old spaces}. To choose these old spaces, the clans play a vital role in several cases.

\vspace{0.1cm}

Suppose we have identified the spaces whose superscripts belong to the set $\{1,...,x\}$. What we have to do now follows from the definition of the configuration.

In configuration 3a, we have $\mu_{1}$ spaces with $s$ new open $k$-sets and the remaining $\mu_{2}$ spaces with $s+1$ new open $k$-sets. Here the spaces with $s$ new open $k$-sets become old spaces: they get as superscripts $\{x+1,...,x+\mu_{1}\}$. (For the moment forget the question \emph{which is which}. In fact, sometimes, the quotient-spaces $Q^{\alpha_{1}}$,  $Q^{\alpha_{2}}$, ...,  $Q^{\alpha_{r}}$ are isomorphic, hence they cannot be distinguished. This happens, for instance, when $\alpha_{1} ^{*}=\alpha_{2} ^{*}=...=\alpha_{r} ^{*}$).

As regards configuration 1a, recall that it is configuration 3a when $\mu _{2}=0$.

\vspace{0.1cm}

In configuration 3d, we have $\mu_{2}$ spaces with $s+1=1$ new open $k$-sets and the remaining $\mu_{1}$ spaces with $s=0$ new open $k$-sets. Here the spaces with $1$ new open $k$-set become old spaces: they will receive as superscripts $\{x+1,...,x+\mu_{2}\}$.

As regards configuration 1d, recall that it is configuration 3d when $\mu _{1}=0$.

\vspace{0.1cm}

In configuration 3c, we have $\mu_{1}$ new spaces with $s$ new open $k$-sets and the remaining $\mu_{2}=n-1-t_{k}-\mu_{1}$ spaces with $s+1$. These $\mu _{2}$ spaces become old. Among those with $s$ new sets, some of them also become old: to choose them, look at the unions of the paired $k-$systems with each lower non-paired $k-$system (more precisely, the unions of open sets yielded by paired $k-$systems with each open set yielded by each lower non-paired $k-$system). The spaces which contain such unions become old; those which do not contain them remain new.

\vspace{0.1cm}

In configuration 3b, we have also $\mu_{1}$ new spaces with $s$ new open $k-$sets and $\mu_{2}=n-1-t_{k}-\mu_{1}$ new spaces with $s+1$. The $\mu_{1}$ spaces become old. Among the $\mu_{2}$ spaces we choose $|R_{i}|+...+|R_{s+1}|$ to become old. For this choice, we notice that they must contain all the upper non-paired $k'-$sets, $k'>k$, which are disjoint from the paired $k-$systems; and they do not contain lower non-paired $k'-$systems, $k'>k$, unless they share points which are new for $k$ with the paired $k-$systems.

Choose also to become old those spaces which contain the open sets yielded by whole clans (or subclans) rooted at a new lower $k-$system. If no distinction is possible, choose among the $\mu_{2}$ spaces the old spaces at will.

\vspace{0.2cm}

Configuration 1b has $s$ paired $k-$systems and it is configuration 3b with $\mu_{1}=0$ or configuration 3c with $\mu_{2}=0$. To choose the $|R_{1}|+...+|R_{s}|$ spaces which become old,  we do as for configuration 3b.

\vspace{0.2cm}

Concerning configuration 1c, we may have no new  $k-$systems or one clean new upper $k-$system with as many points as new spaces. The clean new upper $k-$system yields a new open $k-$set in all old spaces. When it exists, all new spaces become old.

\vspace{0.1cm}

For configurations 2a, 2b, 4a and 4b,  it is easy to identify the $\mu_{1}$ spaces with $s$ new open $k-$sets and the $\mu_{3}$ spaces with $s+2$ new open $k-$sets. All these spaces become \emph{old}. And if, as we did before, $x$ denotes  the number of spaces we have already classified as \emph{old} while checking $k'-$systems with $k'<k$, then we assign to the $\mu_{1}$ new spaces with $s$ new open $k-$sets superscripts from the set $\{x+1, ... , x+\mu_{1}\}$ and to the $\mu_{3}$ new spaces with $s+2$ new open $k-$sets, superscripts from the set $\{x+\mu_{1}+1, ... , x+\mu_{1}+\mu_{3}\}$.

\vspace{0.1cm}

For configuration 4c, sometimes it is difficult to identify, among the spaces with $s+1$ new $k-$systems, those which should receive a superscript from the set $\{x+\mu_{1}+\mu_{3}+1,...,x+\mu_{1}+\mu_{3}+\mu\}$ where $\mu =|R_{i}|+ ... +|R_{s+1}|$. These values are associated to each one of the $\mu$ points in $R_{i}\cup ...\cup R_{s+1}$; remember that we may have $\mu_{2}>\mu $. To find $\mu$, we do as follows: using the notation of \emph{Theorem 4}, remember that $A_{i}\cap ... \cap A_{x}\neq \emptyset $, hence the sets $A_{i}R_{i}$, ..., $A_{x}R_{x}$ belong to a clan which appears in an old space. Let $q$ be the number of new points in the clan, that means, points which appear in open sets with $k$ points, but not in smaller sets. Knowing $q$, we obtain $q-\mu _{3}$ as the number of new points (and also of new spaces) in $A_{i}R_{i}\cup ...\cup A_{s+1}R_{s+1}$.

To identify these $\mu $ spaces among the $\mu_{2}$ new spaces with $s+1$ new open $k-$sets we look at the non-paired $k'-$systems for $k'>k$. Recall that lower systems $\{a_{1}, ..., a_{k'}, n\}$ yield open $k'-$sets in spaces $Q^{j}$ with $j\in \{a_{1}, ..., a_{k'}\}$, and upper systems $\{a_{1}, ..., a_{k'}\}$ yield open $k'-$sets in spaces $Q^{j}$ with $j\notin \{a_{1}, ..., a_{k'}\}$.

Three consequences which help distinguish the $\mu $ spaces among the $\mu_{2}$ ones:

Suppose we have a clan rooted on an upper non-paired system, say $A_{1}R_{1}$: remember Observation 3, and look at upper non-paired systems of this clan; if
$\{a_{1}, ..., a_{j}\}$ is one of them, then it appears in all spaces $Q^{\beta }$ except  in those with $\beta \in  \{a_{1}, ..., a_{j}\}$.

We may also suppose we have clans rooted at upper non-paired $k'-$systems for $k'>k$, yielding clean $k'-$sets and covering all new points for $k'>k$. Such $k'-$sets will all appear in the $\mu $ new spaces we want to distinguish right now but they will never appear together in the remaining new spaces.

Finally, remember Observation 2 and consider now a clan rooted at a lower system, say $A_{i}R_{i}n$ with two descendants, $A_{i}R_{i}X n$ and $A_{i}R_{i}Y n$. With $j\in R_{i}$, we see that $Q^{j}$ contains $A_{i}R_{i}X$ and $A_{i}R_{i}Y$ but for $j\notin A_{i}R_{i}\cup ...\cup A_{s+1}R_{s+1}$ the clan is not there.

\vspace{0.2cm}

Let us illustrate these cases.

\emph{Example 6:} Take $\cal T$ with $T$ having enough points, and the following initial configurations:

\vspace{0.1cm}

\hspace{1cm} $1\ $ ; $23$ ; $14\ $  ; $\ 15$

\hspace{1cm} $1n\ $; $\ \   $  ; $14n\ $; $15n\ $; $16n\ $; $17n\ \ $; ...

\vspace{0.1cm}

For $k=1$ we have configuration 1b.
For $k=2$ we have configuration 4c where spaces $Q^{4}$, $Q^{5}$ and $Q^{\beta }$ for $\beta >7$ cannot be distinguished.

We want to identify spaces $Q^{4}$ and $Q^{5}$ among those new spaces with $s+1$ new open 2-sets.

Suppose we have a clan rooted at $\{2,3\}$ with $\{2,3,8,9\}$ as a descendant. These two sets, $\{2,3\}$ and $\{2,3,8,9\}$ appear in
$Q^{4}$ and $Q^{5}$, but not in $Q^{8}$ or $Q^{9}$, which excludes immediately $Q^{8}$ and $Q^{9}$ from the family of old spaces yielded by the $k-$iteration of our procedure.

As an alternative, suppose we have an upper, non-paired $k'-$system  for $k'> k$, yielding a clean $k'-$set, say $\{8,9,10\}$.
The set $\{8,9,10\}$, even if it is a root of a trivial clan, does not appear in $Q^{8}$, $Q^{9}$ or $Q^{10}$ but it appears in all other spaces, making it possible to distinguish between $Q^{4}$ and $Q^{5}$, on one side, and $Q^{8}$, $Q^{9}$ and $Q^{10}$, on the other side.

Finally, suppose we have a clan rooted at $\{1,5,n\}$ (or at $\{1,n\}$), with $\{1,5,8,n\}$ and $\{1,5,9,n\}$ as descendants. These both sets appear in $Q^{5}$ (together with $\{1,5\}$ and $\{1\}$) but they do not appear together neither in $Q^{8}$ nor $Q^{9}$, which allows us to recognize $Q^{5}$ as one of the $\mu$ spaces for $k=2$.

\vspace{0.1cm}

If we have no way to recognize the $\mu$ spaces which become old among those with the $s+1$ new open $k-$sets, then we are free to choose them as we wish. In our example, if nothing distinguishes $Q^{4}$, $Q^{5}$ and $Q^{8}, ...$, then we choose two of them to be $Q^{4}$ and $Q^{5}$, and to become old. The other ones are new spaces for $k=3$ and later they'll become $Q^{8}, ..., Q^{n-1}$.

\vspace{0.2cm}

Just a remark concerning the question \emph{which is which} we mentioned above. This question is irrelevant for the reconstruction. However, if we are interested, an analysis of the clans allows us to choose the superscript of each space among those previously assigned as possible. With the help of the clans, we may also assign the right number to the symbols which represent the points of each quotient space.

\section{A few open questions}

It would be desirable to know a necessary and sufficient condition for the existence of a solution to the \emph{Problem} stated in the Introduction. We know of no such condition but, meanwhile, look at the following preliminary facts:

1. Given a topological space $(T, \cal T)$, the $n-1$ spaces formed by topological identification of $n$ with each one of the other points yield, for $k=1,...,n-1$, one of the 13 configurations we have listed and the successive configurations never infringe the Observations 1 through 5 about clans.

2. Reciprocally, for the \emph{Problem} to have at least one solution, the following is a minimal set of conditions that must be satisfied: For $k=1,...,n-1$, the new covering sets (or the unions of a new covering set with $\{n\}$) which we successively obtain with our algorithmic procedure always form one of the 13 configurations that we have listed in Section 4, and the successive configurations never infringe the Observations 1 through 5 about clans.

\vspace{0.1cm}

\textbf{A conjecture:} In this paper, we can verify that the only cases where more than one solution was reached, were cases where trivial clans rooted at some $k-$set with $k<n-1$ exist (including \emph{Examples 1 and 2}). We conjecture that to have more than one solution (up to homeomorphism) we must have trivial clans rooted at $k-$sets for some $k<n-1$.

\vspace{0.2cm}

\textbf{Enumerative questions} may also be asked. For instance:

1. How many topologies, distinct up to homeomorphism, can we define on $T=\{1, ..., n\}$ with $n^{*}=T$?
2. And how many have $n^{*}=\{n\}$?

These two particular cases of topologies can be considered extreme cases for the covering set $n^{*}$ of point $n$. For them, it is an immediate consequence of the definitions that the following two statements are valid:

1. When there is a solution where $n^{*}=T$, only 1a, 1c and 3a may appear as configurations, and also  1b  but just for $k=n-1$;

2. When there is a solution where $n^{*}=\{n\}$, only 1b, 1c, 1d, 3c and 3d may appear as configurations.

\section{Some more examples}

Recall \emph{Example 3} at the end of Section 4 and the topologies of its three quotient-spaces. For $k=1$, we have $\mu_{1}=2$ spaces with $s=1$ system, and $\mu_{2}=1$ space with $s+1=2$ systems and $p=2$ new points. This is a configuration of type 3. Since $\mu_{1}=p$, configuration 3b is excluded. We distinguish 3a from 3c as we explained in Section 6: we choose configuration 3a, that means sets $\{1\}$ and $\{2\}$ as 1-systems. Let $Q^{1}$ and $Q^{2}$ be the spaces with topologies $U_{1}$ and $U_{2}$ respectively. $U_{3}$ is the topology of the new space; here, for $k=2$, $\{g,h\}$ is an old doublet, hence we have $\mu_{1}=1$ space with $s=0$ new doublets and $\mu_{2}=\mu_{3}=0$. It is configuration 1c. For $k=3$ we have $\mu_{1}=1$ space with $1$ new mixed triplet $A_{1}R_{1}$ where $A_{1}=\{g,h\}$ and $R_{1}=\{j\}$.
\vspace{0.1cm}

\emph{Example 7:} Take a topology on $T=\{1,2,3,4\}$, with $\{1\}$ and $\{2\}$ as singletons and $\{2,3\}$ as a new $2-$system. Besides $\{1,2,3\}$, the space $Q^{1}$ will have as open sets $\{2\}$ and $\{2,3\}$; $Q^{2}$ will have as open set only $\{1\}$; and $Q^{3}$ will have as open sets $\{1\}$, $\{2\}$ and $\{1,2\}$. In symbols, as given by the data of the problem, we may be told that one of the spaces (later it will be named $Q^{1}$) has $a, ab, abc$ as covering sets, another space (it will be named $Q^{2}$) has $f, fgh$ as covering sets, a third space (to be named $Q^{3}$) has $y, z, yz, yzw$ as open sets ($yz$ is an old $2-$system). For $k=1$, we have $\mu_{1}=2$, $\mu_{2}=1$ and configuration 3a, as explained in Section 6; for $k=2$, configuration 1a. While looking for new $2-$systems in the new spaces (here $Q^{3}$ is the only new space for $k=2$) we don't find the new $2-$system which is $\{2,3\}=A_{1}R_{1}$; it appears however in the old space $Q^{1}$ as $ab$. In this case, $\{1,2,3\}$ is an old $3-$system.

\vspace{0.1cm}

Now, one case of configuration 1a, with $s>1$, where no point $x$ exists such that  $x\notin A_{1}R_{1}\cup ... \cup A_{s+1}R_{s+1}$. (Here $A_{1}=...=A_{s+1}=\emptyset $.)

\vspace{0.1cm}

\emph{Example 8:} Take ${\cal T'}$ defined on $T=\{1,...,7\}$ with covering sets $\{1,2\}$, $\{3,4\}$, $\{5,6\}$ and $T=\{1,...,7\}$. We have configuration 1c for $k=1$, and 1a for $k=2$. Compare with ${\cal T''}$ with two paired open sets $\{1,2\}$, $\{1,2,7\}$, $\{3,4\}$, $\{3,4,7\}$ and $T=\{1,...,7\}$. We have here configuration 1c for $k=1$ and 1b for $k=2$. The quotient-spaces in both cases can be written $Q^{1}=Q^{2}=\{ab, cd, abcd, abcdef\}$,
$Q^{3}=Q^{4}=\{ab, cd, abcd, abcdef\}$, $Q^{5}=Q^{6}=\{ab, cd, abcd, abcdef\}$. Obviously, the names of these symbols will be different:

For ${\cal T'}$, $Q^{1}=Q^{2}=\{34, 56, 3456, 345612\}$, $Q^{3}=Q^{4}=\{12, 56, 1256, 125634\}$, $Q^{5}=Q^{6}=\{12, 34, 1234, 123456\}$.

For ${\cal T''}$, we have $Q^{1}=...=Q^{6}=\{12, 34, 1234, 123456\}$.

Note that, for $k=2$, if we choose configuration 1a, then, in ${\cal T'}$, $n^{*}=T$; if we choose configuration 1b, then, in ${\cal T''}$, $n^{*}=\{n\}=\{7\}$. Note also that the clans of covering sets in these quotient-spaces are all trivial.

\end{document}